\documentclass[reqno]{amsart}
\theoremstyle{plain}
\newtheorem{thm}{\it Theorem}[section]
\newtheorem{prop}[thm]{\it Proposition}

\newtheorem{lem}[thm]{\it Lemma}
\theoremstyle{remark}
\newtheorem{defn}[thm]{Def{}inition}
\newtheorem{rem}[thm]{Remark}
\newtheorem{exa}[thm]{Example}
\numberwithin{equation}{section}
\usepackage{enumerate}

\newcommand{\NN}{\mathbb{N}}

\def\cc#1{\{#1\}}
\def\pp#1{\|#1\|}
\begin{document}
\baselineskip08pt
\title{On $\mathfrak{I}$-Reconstruction Property}
\author[L. K. Vashisht and Geetika Khattar]{ L. K. Vashisht and Geetika Khattar\\D\MakeLowercase{epartment} \MakeLowercase{of}
M\MakeLowercase{athematics},\\ U\MakeLowercase{niversity of} D\MakeLowercase{elhi}, D\MakeLowercase{elhi 110007}, I\MakeLowercase{ndia}\\
E\MakeLowercase{-mail: lalitkvashisht@gmail.com ,
geetika1684@yahoo.co.in}}
\begin{abstract}\baselineskip8pt
 Reconstruction
property in Banach spaces introduced and \break  studied by Casazza
and Christensen in [1]. In this paper we introduce \break
reconstruction
 property in Banach spaces which satisfy  $\mathfrak{I}$-property.
 A \break characterization of  reconstruction property in Banach spaces which satisfy  \break $\mathfrak{I}$-property in
 terms of frames in Banach spaces is obtained.  Banach frames associated with reconstruction property
  are discussed.
 \vskip.25em\noindent


Reconstruction property.
\end{abstract}

\maketitle \thispagestyle{empty} \baselineskip12pt

\section{Introduction}
Let $\mathcal{H}$ be an infinite dimensional separable complex
Hilbert space  with inner product $\langle ,     \rangle$. A system
$\{f_k\} \subset \mathcal{H}$ called a \emph{frame}(Hilbert) for
$\mathcal{H}$ if there exists positive constants $\mathrm{A}$ and
$\mathrm{B}$  such that
$$\mathrm{A}\|f\|^2\leq \sum_{k=1}^{\infty} |\langle f, f_k\rangle |^2 \leq
\mathrm{B}\|f\|^2, \  f \in \mathcal{H}. $$\

The positive constants $\mathrm{A}$ and $\mathrm{B}$ are called
\emph{lower} and \emph{upper bounds} of the frame $\{f_k\}$,
respectively. They are not unique.

The operator  $T:\ell^2\rightarrow \mathcal{H}$ given by
$T(\{c_k\})=\sum_{k=1}^{\infty}c_kf_k, \{c_k\} \in \ell^2$ is called
the \emph{synthesis operator} or \emph{pre-frame operator}. Adjoint
of $T$ is given by \break $T^*:\mathcal{H} \rightarrow \ell^2$,
$T^*(f)= \{ \langle f, f_k \rangle\} $ and is called the \emph{analysis
operator}. Composing $T$ and $T^*$ we obtain the \emph{frame
operator} $S=TT^{*}:\mathcal{H}\rightarrow \mathcal{H}$ given by
\break  $S(f)=\sum_{k=1}^{\infty}\langle f, x_k\rangle f_k, f \in
\mathcal{H}$. The frame operator $S$ is a positive continuous
invertible linear operator from $\mathcal{H}$ onto $\mathcal{H}$.
Every vector $f \in \mathcal{H}$ can be written as:
$$f = SS^{-1}f =\sum_{k=1}^{\infty}< S^{-1}f, f_k > f_k.\ (Reconstruction \ formula) $$
The series in the right hand side converge unconditionally and is
called \emph{reconstruction formula} for $\mathcal{H}$.  The
representation of $f$ in reconstruction formula need not be unique.
Thus,  frames are \emph{redundant systems} in a Hilbert space which
yield one natural representation for every vector in the concern
Hilbert space, but which may have infinitely  many different
representations for a given vector.\

Duffin and Schaeffer in [2] while working in non-harmonic Fourier
series developed an abstract framework for the idea  of
time-frequency atomic decomposition by Gabor[3] and defined
\emph{frames} for Hilbert spaces. Due to some reason the theory of
frames was not continued until 1986 when the fundamental work of
Daubechies, Grossmann and Meyer  published in [4].
Gr$\ddot{o}$chenig in [5] generalized Hilbert frames to Banach
spaces. Before the concept of Banach frames was formalized, it
appeared in the foundational work of Feichtinger and
Gr$\ddot{o}$chenig [6, 7] related to \emph{atomic decompositions}.
Atomic decompositions appeared in the field of applied mathematics
providing many applications [8, 9]. An atomic decomposition allow a
representation of every vector of the space via a series expansion
in terms of a fixed sequence of vectors which we call \emph{atoms.}
On the other hand Banach frame for a Banach space ensure
reconstruction via a bounded linear operator or \emph{synthesis
operator}. Frames play an important role in the theory of nonuniform
sampling[10], wavelet theory [11, 12], signal processing [2, 10],
and many more. For a nice introduction of frames and their technical
details one may refer to [13].\

During the development of frames and expansions systems in Banach
spaces Casazza and Christensen introduced  \emph{reconstruction
property} for Banach spaces in [1]. Reconstruction property is an
important tool in several areas of mathematics and engineering.
Infact, it is related to bounded approximation property. Casazza and
Christensen in [1] study perturbation theory related to
reconstruction property.  They develop more general perturbation
theory that does not force equivalence of the sequences.\

In this paper we introduce and study reconstruction
 property in Banach spaces which satisfy  $\mathfrak{I}$-property.
 A characterization of  $\mathfrak{I}$-reconstruction property in terms of frames in Banach spaces is obtained.
  Banach frames associated with reconstruction property
  are discussed.

\section{Preliminaries}
Throughout this paper $\mathcal{X}$ will denotes an infinite
dimensional Banach space over
 a field $\mathbb{K}$ (which can be $\mathbb{R}$ or $\mathbb{C}$ ),
 $\mathcal{X}^*$ be the conjugate space, and for a sequence $\{f_k\}\subset \mathcal{X}$,
 $[f_k]$ denotes closure of $spanf_k$ in norm topology of
 $\mathcal{X}$. The map $\pi:\mathcal{X} \rightarrow
 \mathcal{X}^{**}$ denotes the canonical mapping from $\mathcal{X}$ into
 $\mathcal{X}^{**}$.

\begin{defn}[\cite{5}]
 Let $\cc{f^*_k}\subset
\mathcal{X}^*$ and $\mathfrak{S}:\mathcal{X}_d\to \mathcal{X}$ be
given, where $\mathcal{X}_d$ is an associated Banach space of scalar
valued sequences. A system  $\mathcal{G}
\equiv(\cc{f^*_k},\mathfrak{S})$ is called a \emph{Banach frame} for
$\mathcal{X}$ with respect to $\mathcal{X}_d$ if
\begin{enumerate}[(i)]
\item $\cc{f^*_k(f)}\in \mathcal{X}_d$, for each $f\in \mathcal{X}$.
\item There exist positive constants $\mathrm{C}$ and $\mathrm{D}$  $(0<\mathrm{C} \le \mathrm{D} <\infty)$ such that
\begin{align}
\mathrm{C}\pp{f}_\mathcal{X}\le
\pp{\cc{f^*_k(f)}}_{\mathcal{X}_d}\le
\mathrm{D}\pp{f}_{\mathcal{X}},  \  for \   each \ f\in \mathcal{X}.
\end{align}
\item
$\mathfrak{S}$ is a bounded linear operator such that
\begin{align*}
\mathfrak{S}(\cc{f^*_k
(f)})=f, \  for \ each \ f\in \mathcal{X}\,.
\end{align*}
\end{enumerate}
\end{defn}
As in case of frames for a Hilbert space, positive constants
$\mathrm{C}$ and $\mathrm{D}$ are called {\it lower} and {\it upper
frame bounds} of the Banach frame $\mathcal{G}$, respectively. The
operator $\mathfrak{S}:\mathcal{X}_d\to \mathcal{X}$ is called the
{\it reconstruction operator} (or the {\it pre-frame operator}). The
inequality $2.1$ is called the {\it frame inequality}.

The Banach frame $\mathcal{G}$ is called {\it tight} if
$\mathrm{C}=\mathrm{D}$ and {\it normalized tight} if
$\mathrm{C}=\mathrm{D}=1$. If there exists no reconstruction
operator $\mathfrak{S}_m$ such that $(\{f^*_k\}, \mathfrak{S}_m)_{k
\ne m} \break (m \in \mathbb{N})$ is  Banach frame for
$\mathcal{X}$, then $\mathcal{G}$ will be called an \emph{exact}
Banach frame.

The notion of retro Banach frames introduced and studied in [14].

\begin{defn}[\cite{14}]
A system $\mathcal{F} \equiv (\{f_k\},\Theta)$ $(\{f_k\}\subset
\mathcal{X },\Theta:\mathcal{Z}_d \rightarrow X^*)$ is called a
\emph{retro Banach frame} for $\mathcal{X}^*$ with respect to an
associated sequence space $\mathcal{Z}_d$ if
\begin{enumerate}[(i)]
\item  $\{f^*(f_k)\}\in \mathcal{Z}_d$, for each $f^*\in \mathcal{X}^*$.
\vskip 2pt
\item There exist positive constants $(0< \mathrm{A}_0 \le \mathrm{B_0}< \infty)$ such that
$$\mathrm{A_0}\|f^*\| \leq \|\{f^*(f_k)\}\|_{\mathcal{Z}_d} \leq\mathrm{B_0}\|f^*\|, \ for \ each \ f^* \in\mathcal{X}^*.$$
\vskip 2pt
\item $\Theta:$ $\{f^*(f_k)\} \rightarrow$ $f^*$ is a bounded linear operator from  $\mathcal{Z}_d$ onto $\mathcal{X}^*.$
\end{enumerate}

\end{defn}
The positive constant $\mathrm{A}_0,  \mathrm{B_0}$ are called
\emph{retro frame  bounds} of $\mathcal{F}$ and
operator $\Theta:\mathcal{Z}_d \rightarrow X^*$ is called
\emph{retro pre-frame operator} (or simply \emph{reconstruction operator})
associated with $\mathcal{F}$.

\begin{lem}.
Let $\mathcal{X}$ be a Banach space and $\left\{f_n^*\right\}\subset
\mathcal{X}^*$ be a sequence such that $\left\{f \in
\mathcal{X}:f_n^*\left(f\right) = 0, \text{for all}\; n\in
\mathbb{N}\right\} = \left\{0\right\}$. Then, $\mathcal{X}$ is
linearly isometric to the \break Banach space $\mathcal{Z} =
\left\{\left\{f_n^*\left(f\right)\right\}:f\in \mathcal{X}\right\}$,
where the norm  is  given by
$\left\|\left\{f_n^*(f)\right\}\right\|_\mathcal{Z}=\left\|f\right\|_\mathcal{X}
,f\in \mathcal{X}$.
\end{lem}
Casazza and Christensen in [1] introduced reconstruction property in
Banach spaces.
\begin{defn}[\cite{1}] Let $\mathcal{X}$ be a separable Banach space. We say that a
sequence $\{f_k^*\} \subset \mathcal{X}^*$ has the
\emph{reconstruction property} for $\mathcal{X}$ with respect  to
$\{f_k\}\subset \mathcal{X}$, if
$$ f =  \sum_{k=1}^{\infty}f_k^*(f)f_k \ \  \ for \ all \  f \in \mathcal{X}. \ \ \ \ \ $$
\end{defn}

In short, we will also say $(\{f_k\}, \{f_k^*\})$ has reconstruction
property for $\mathcal{X}$. More precisely, we say that  $(\{f_k\},
\{f_k^*\})$ is a \emph{reconstruction system} for $\mathcal{X}$.

\begin{rem}
A interesting example for a reconstruction property is given in [1]:
Let $\{f_k^*\} \subset \ell^{\infty}$ and  $\{f_k^*\}$ is unitarily
equivalent to the unit vector basis of $\ell^2$. Then,  $\{f_k^*\}$
has a reconstruction property with respect to its own pre-dual (that
is, expansions with respect to the orthonormal basis). Further
examples on  reconstruction property are discussed in  Example
$3.4$.
\end{rem}
\begin{defn}
A reconstruction system $(\{f_k\}, \{f_k^* \}) $  for $\mathcal{X}$
is said to be
\begin{enumerate}[(i):]
\item \emph{pre-shrinking} if $[f^*_k]= \mathcal{X}^*$.
\item \emph{shrinking} if $(\{f^*_k\}, \{f_k \}) $ is a
reconstruction system for $\mathcal{X^*}$.

\end{enumerate}

\end{defn}

Regarding existence of Banach spaces which have  reconstruction
system, Casazza and Christensen proved the following result.

\begin{prop}[\cite{1}]
There exists a Banach space $\mathcal{X}$ with the following
properties:
\begin{enumerate}[$(i)$:]
\item There is a sequence $\{f_k\}$ such that each $f\in \mathcal{X}$
has a expansion $f =  \sum_{k=1}^{\infty}f_k^*(f)f_k. $
\item $\mathcal{X}$ does not have the reconstruction property with
respect to any pair $(\{h_k\}, \{h_k^*\}).$

\end{enumerate}
\end{prop}

The notion of reconstruction property is related to Bounded
Approximation Property(BAP). If $(\{f_k\}, \{f_k^*\})$ has
reconstruction property for $\mathcal{X}$, then $\mathcal{X}$ has
the bounded approximation property. So, $\mathcal{X}$ is isomorphic
to a complemented subspace of a Banach space with a basis. It is
also used to study geometry of Banach spaces. For more results and
basics on reconstruction property and bounded approximation property
one may refer to [15] and references therein.

 \section{$\mathfrak{I}$ - Reconstruction Property}
\begin{defn} Suppose $\{f_k^*\} \subset \mathcal{X}^*$ has the
\emph{reconstruction property} for $\mathcal{X}$ with \break respect
to $\{f_k\}\subset \mathcal{X}$. Then, we say that  $(\{f_k\},
\{f_k^* \}) $ satisfy \emph{property $\mathfrak{I}$}
  if \break $\inf\limits_{1 \leq k < \infty} \|f_k\|> 0$
 and there exists a  functional $\Psi$ $\in$ $\mathcal{X}^{**}$ such that
$\Psi(f_k^*) = 1$, for all $k \in \mathbb{N}$. In this case we  say
that  $(\{f_k\}, \{f_k^* \}) $ is a \emph{$\mathfrak{I}$ -
reconstruction system} for $\mathcal{X}$.

\end{defn}

\begin{rem}
If $\inf\limits_{1 \leq k < \infty} \|f_k\| = 0$ and there exists a
functional $\Psi$ $\in$ $\mathcal{X}^{**}$ such that $\Psi(f_k^*) =
1$, for all $k \in \mathbb{N}$,then we say that $(\{f_k\}, \{f_k^*
\}) $ is a  \emph{$\mathfrak{I_\omega}$-reconstruction system}
(or\emph{ weak $\mathfrak{I}$ - reconstruction system} for
$\mathcal{X}$).
\end{rem}
\begin{rem}
A $\mathfrak{I}$-reconstruction system is actually a dual system of
a $\Phi$-Schauder frame $[16]$ in the context of reconstruction
property.
\end{rem}

\begin{exa} Let $\mathcal{X} = c_0$ and
$\{e_k\}\subset \mathcal{X}$
 be a sequence of canonical unit vectors. Define $\{f_k^*\} \subset
 \mathcal{X}^*$ by
$f_1^*(f) = \frac{1}{2}\xi_{1}, \ f_2^*(f)=\frac{1}{2}\xi_{1}, \
f_k^*(f) = \xi_{k-1}, \ \ f=\{\xi_k\} \in \mathcal{X}$. Then,
$\{f_k^*\}$ has a reconstruction property with respect to $\{f_k\}
\subset \mathcal{X}$, where \break $f_1 = e_1, \ f_2 = e_1,  f_k =
e_{k-1}.$  Hence $(\{f_k\}, \{f_k^* \})$ is a $\mathfrak{I}$ -
reconstruction system for $\mathcal{X}$ [See Proposition $3.5$].
Note that the reconstruction system $(\{f_k\},
\{f_k^* \})$ is shrinking.\\
 Now define  $\{g_k^*\} \subset
 \mathcal{X}^*$ by
$g_1^*(f) = \xi_1,  \ g_k^*(f) = \xi_{k-1}, f=\{\xi_k\} \in
\mathcal{X}$. Then, $\{g_k^*\}$ has a reconstruction property with
respect to $\{g_k\} \subset \mathcal{X}$, where $g_1 = 0, \ g_k =
e_{k-1}.$ By Proposition $3.5,$ $(\{g_k\}, \{g_k^* \})$ is not a
$\mathfrak{I}$ - reconstruction system for $\mathcal{X}$. Note that
$(\{g_k\}, \{g_k^* \})$ is $\mathfrak{I}_\omega$-reconstruction
system which is shrinking.
 Thus, a shrinking reconstruction system for $\mathcal{X}$ need not be a $\mathfrak{I}$-reconstruction system.

\end{exa}

We now give   a characterization of  a $\mathfrak{I}$ -
reconstruction system for $\mathcal{X}$ as claimed in section $1,$
in terms of frames.
\begin{prop}
Let $(\{f_k\},\{f_k^*\})$ be a reconstruction system for
$\mathcal{X}$ with \break $\inf\limits_{1\leq k\leq \infty}
\parallel f_k\parallel \
> 0 $. Then , $(\{f_k\}, \{f_k^*\})$ satisfy property $\mathfrak{I}$ if and only if  there is no
retro pre-frame operator $\Theta_{0}$ such that  $( \{f_k^*-
f_{k+1}^*\},\Theta_0)$ is retro Banach frame for $[f_k^*]^*$.
\end{prop}
This is an immediate consequence of the following lemma.
\begin{lem}
Let $(\{f_k\},\{f_k^*\})$ be a pre-shrinking reconstruction system
for $\mathcal{X}$. Then, $( \{f_k\}, \{f_k^*\})$ is a
$\mathfrak{I_\omega}$-reconstruction system if and only if there
exists no retro pre-frame operator $\widehat{\Theta}$ such that
 $( \{f_k^*- f_{k+1}^*\},\widehat{\Theta})$ is retro Banach frame for $\mathcal{X}^{**}$.\
\end{lem}
\proof Forward part is obvious. Indeed, by using lower retro frame
inequality of  $( \{f_k^*- f_{k+1}^*\},\widehat{\Theta})$ and
existence of $\Psi \in \mathcal{X}^{**}$ such that $\Psi(f_k^*-
f_{k+1}^*) = 0,$ for all  $k \in \mathbb{N},$ we obtain $\Psi = 0.$
This is a
contradiction.\\
For reverse part, let if possible,  there is no reconstruction
operator $\widehat{\Theta}$ such that \break $( \{f_k^*-
f_{k+1}^*\},\widehat{\Theta})$ is a retro Banach frame for $\mathcal{X}^{**}$.
Then, Hahn Banach \break Theorem
 force to admit a non zero functional $\phi\in E^{**}$ such that
$\phi(f_k^* - f_{k+1}^*) = 0$, for all $k \in \mathbb{N}$. That is,
$\phi(f_k^*) =\phi(f_{k+1}^*)$, for all $k \in \mathbb{N}$. Put
$\phi(f_k^*)=\alpha$,  for all $k \in \mathbb{N}$. If $\alpha= 0,$
then $\phi(f_k^*)=0$ for all $\ k \in \mathbb{N}.$ But $(\{f_k\},
\{f_k^*\})$ is pre-shrinking, therefore $\phi$ = 0 , a
contradiction. Thus $\alpha$ $\neq$  0. Put $\Psi$=
$\frac{\phi}{\alpha}$. Then $\Psi$ $\in$ $\mathcal{X}^{**}$ is such
that $\Psi(f_k^*) = 1$   for all $ \ k \in \mathbb{N}.$ Thus,
$(\{f_k\}, \{f_k^*\})$ is a $\mathfrak{I_\omega}$- reconstruction
system.
\endproof

\begin{rem}
Note that Lemma $3.6$ is no longer true if $(\{f_k\},\{f_k^*\})$ is
not pre-shrinking.
\end{rem}

\emph{Application}: Let $\mathcal{X}= L^2(a, b)$. Consider a
boundary value problem(BVP) with a set of n boundary conditions:\

BVP:$   \ \ \nabla(f)= \lambda f , \ \ \ \Xi(f)=0,$\\
where  $\nabla(\bullet)= (\bullet)^n + \Phi_1(\xi)(\bullet)^{n-1} +
..... + \Phi_n(\xi)(\bullet)$  is a linear differential operator
with
$\Phi_j \in C^{n-k}[a, b],$ and $\Xi(f)=0$ denotes the set of n boundary conditions: $\Xi_j(f) = \Sigma_{k=1}^{n} [\alpha_{j, k}\Phi^{k-1}(a) + \beta_{j, k}\Phi^{k-1}(b)] = 0.$\\
It is given in [17] (at page 66) that for a large class of boundary
conditions (which are known as regular
 boundary conditions), the BVP  admits a system $\{\Phi_n(\xi)\} $ and $\{\Psi_n(\xi)\}$
 consisting of eigenfunction associated with given BVP such that\

$\Phi_n(\xi) = A_n [cos\frac{2\pi n \xi}{b-a} + O(\frac{1}{n})]; \ \ $\\

  $\Psi_n(\xi) = B_n [sin\frac{2\pi n \xi}{b-a} + O(\frac{1}{n})], n = 0,1,2,3,....$\\
It is well known that the corresponding to $\{f_n\} \equiv
\{cos\frac{2\pi n \xi}{b-a}\} \bigsqcup  \{sin\frac{2\pi n
\xi}{b-a}\}$ there  exists a $\{f^*_n\} \in \mathcal{X}^*$ such that
$(\{f_n\}, \{f^*_n\})$ is a reconstruction system for  $\mathcal{X}
= L^2(a, b).$
 Now
 $\|\Phi_n(\xi) - A_n cos\frac{2\pi nt}{b-a}\|^2 = O(\frac{1}{n^2}) $
 and
  $\|\Psi_n(\xi) - B_n sin\frac{2\pi nt}{b-a}\|^2 = O(\frac{1}{n^2}). $\\
 Therefore, by using  Paley and Wiener theorem in [18 at page 208], there exists a sequence $\{f^*_n\} \in
 \mathcal{X}^*$ such that $\{f^*_n\}$ admits a reconstruction system
 with respect to $\{\Phi_n\} \bigsqcup \{\Psi_n\}$. This
 reconstruction system is not of type $\mathfrak{I}_\omega$.
 Therefore, by using Lemma $3.6$, there exists a retro pre-frame operator $\widehat{\Theta}$ such that
 $( \{f_n^*- f_{n+1}^*\},\widehat{\Theta})$ is retro Banach frame for $\mathcal{X}^{**}$.
Recall that if we write a function in terms of reconstruction
system, then computation of all the coefficients is required. If
calculation of coefficients which appear in the series expansion of
a given reconstruction system are complicated, then we reconstruct
the function by
pre-frame operator of  $( \{f_n^*- f_{n+1}^*\},\widehat{\Theta})$.\\

 The following proposition provides a sufficient condition for
 a reconstruction system to satisfy property $\mathfrak{I_\omega}$.
\begin{prop} Let $(\{f_k\}, \{f_k^*\})$ be a  reconstruction system for $\mathcal{X}$.
 If there exists a vector $f_0$ in $\mathcal{X}$ such
that ${f_k^*(f_0)} = 1$ for all $ \ k \in \mathbb{N}$, then
$(\{f_k\},\{f_k^*\})$ is a \break $\mathfrak{I_\omega}$-
reconstruction system.
\end{prop}

\proof Let $\pi : \mathcal{X} \ \longrightarrow \mathcal{X}^{**}$ be
the canonical embedding of $\mathcal{X}$ into $\mathcal{X}^{**}$.
Then $\Psi=\pi(f_0)\in \mathcal{X}^{**}$ is such that
$\Psi(f_k^{*})=1$, for all $ \ k \in\mathbb{N}$. Thus, $(\{f_k\},
\{f_k^*\})$  is a $\mathfrak{I_\omega}$-reconstruction system for
$\mathcal{X}$.
\endproof

\begin{rem}
The condition in Proposition $3.8$ is not necessary. However, if
$\mathcal{X}$ is reflexive, then the condition given in Proposition
$3.8$ turns out to be necessary. Moreover, this is equivalent to the
condition: \emph{There exists no pre-frame operator $\mathfrak{S_0}$
such that $(
\{f_k^*- f_{k+1}^*\},\mathfrak{S_0})$ is a Banach frame for $\mathcal{X}$}. \\
\end{rem}

To conclude the section we show that a given
$\mathfrak{I}$-reconstruction system  in \break  Banach spaces
produce another $\mathfrak{I}$-reconstruction system:
 Consider a $\mathfrak{I}$-reconstruction  system $(\{f_k\}, \{f_k^*\})$ for $\mathcal{X}$.\\
Let $\mathcal{U} = \{\{\gamma_i\} \subset \mathbb{K}:
\sum_{i=1}^{\infty}{\gamma_{i}f_i} \ converges\}.$
\\
Then  $\mathcal{U}$ is a Banach space with norm given by
\\
$$\parallel\{\gamma_{i}\}\parallel_{\mathcal{U}}= \sup_{1\leq k\leq \infty}
\parallel\sum_{i=0}^{k}{\gamma_{i}f_i}\parallel_{\mathcal{X}}.$$
\\
\vskip.25em\noindent Define  $\Gamma : \mathcal{X} \longrightarrow
\mathcal{U}$ by $\Gamma(f)$ = $\{f_k^*(f)\} , f \in \mathcal{X}$.
\\
Then $\Gamma$ is an isomorphism of $\mathcal{X}$ into $\mathcal{U}.$
\\ Also $\widetilde{\Theta} : {\mathcal{U}}\longrightarrow \mathcal{ X}$ defined by
$\widetilde{\Theta}(\{\gamma_i\})$ =$
\sum_{i=1}^{\infty}{\gamma_{i}f_i}$ is also a bounded linear
operator from $\mathcal{U}$ onto $\mathcal{X}$.
\\
Put $\Xi=Ker \widetilde{\Theta}$. Then $\Xi$ is a closed subspace of
$\mathcal{U}$ such that  $\Gamma(\mathcal{X})\bigcap \Xi= \{0\}.$
Moreover, if  $\{\gamma_i\}\in \mathcal{U}$ is any element such
 that $f \ = \sum_{i=1}^{\infty}{\gamma_{i}f_i}$, then $\{f_i^*(f)\}\in\ \Gamma(\mathcal{X})$ and
\\
$$\sum_{i=1}^{\infty}(\gamma_{i}-\{f_i^*(f)\}) \ {f_i}
=\sum_{i=1}^{\infty}{\gamma_{i}f_i}  -  \sum_{i=1}^{\infty}{f_i(f)f_i}$$
$$= 0.$$\\
 Therefore, $(\gamma_{i}-\{f_i^*(f)\})\in \Xi$  is such that $\{\gamma_i\} = \{f_i(f)\}+\{\gamma_i - f_i(f)\}.$
\\
Hence $\mathcal{U}$ = $\Gamma(\mathcal{X})\bigoplus \Xi.$\\
 Let V be projection on $\mathcal{U}$  onto $\Gamma(\mathcal{X})$.\\
 Then, $ V (\{\gamma_i\}) =
\{{f^*_k(\sum_{i=1}^{\infty}\gamma_{i}f_i)}\}, \  \{\gamma_i\}\in
\mathcal{U}$. Therefore, for each k$\in \mathbb{N}$, we have
\\
$$ V(e_k)= \{f^*_k(\sum_{i=1}^{\infty}{\delta_{i, k}f_i)}\} , \ where  \ \delta_{i,k}=(1 \ \  i=k  \ and \  0 \ \ i\neq k)$$
That is: $V(e_k)$ =  $\Gamma(f_k)$ for all $k.$ So, $f_k =
\Gamma^{-1}(V(e_k))$ for all $\  k \  \in \  \mathbb{N}$, where
$\{e_k\}$ is sequence of canonical unit vectors in  $\mathcal{U}$.
 Hence $( \{\Gamma^{-1}(V(e_k)\},\{f^*_k\})$ is a reconstruction system for $\mathcal{X}$ which satisfy property
 $\mathfrak{I}$.\\
 This is summarized in the following proposition.
\begin{prop}
Let $(\{f_k\},\{f_k^*\})$ be a $\mathfrak{I}$-reconstruction system
for $\mathcal{X}$. Then, there exists $(\Gamma^{-1}(V(e_k)) \subset
\mathcal{X}$ such that  $( \{\Gamma^{-1}(V(e_k))\}, \{f^*_k\})$ is a
$\mathfrak{I}$-reconstruction \break system for $\mathcal{X}$, where
$\Gamma$ and $V$are same as in above discussion.

\end{prop}
\section{Associated Banach Frames }

\begin{defn}
Suppose that $\{f_k^*\}$ has the reconstruction property for
$\mathcal{X}$ with \break respect to $\{f_k\}\subset \mathcal{X}$.
Then, there exists a reconstruction operator  $\mathfrak{S}:
\mathcal{X}_d \rightarrow \mathcal{X}$ such that $(\{f_k^*\},
\mathfrak{S})$ is a Banach frame for $\mathcal{X}$ with  respect to
some $\mathcal{X}_d$. We say that $(\{f_k^*\}, \mathfrak{S})$ is an
\emph{associated Banach frame}
 of  $(\{f_k\} ,\{f^*_k\})$.

\end{defn}

 Consider a reconstruction system $(\{f_k\},
\{f^*_k\})$ for a  Banach space $\mathcal{X}$. We can write each
element of $\mathcal{X}$(we can reconstruct $\mathcal{X}$) by mean
of an infinite series formed by $\{f_k\})$ over scalars
$(\{f^*_k(f)\}$. For a
non zero functional $h^*$(say), in general, there is \\
$\bullet$ no $\{h_k\} \subset \mathcal{X}$ such that $(\{f^*_k+
h^*\}$ has the reconstruction  property for $\mathcal{X}$ with
respect to $\{h_k\}$.\\
$\bullet$ no reconstruction operator $\mathfrak{S}_0$ such that
$(\cc{f^*_k+h^*}, \mathfrak{S}_0)$ is a Banach frame for $\mathcal{X}$.\\

More precisely, two natural and important problem arise, namely,
\break existence of
   $\{h_k\} \subset \mathcal{X}$ such that $\{f^*_k+
h^*\}$ has the reconstruction  property for $\mathcal{X}$ with
respect to $\{h_k\}$ and other is the existence of a reconstruction
operator $\mathfrak{S}_0$ associated with $\{f^*_k+ h^*\}$. Cassaza
and Christensen in [1] study some stability of reconstruction
property in Banach spaces in terms of closeness of certain sequence
to a given reconstruction system. In the present section we focus on
pre-frame operator associated with
$\cc{f^*_k+h^*}$.\\
\emph{Motivation}:  Consider a signal space $\mathcal{H}_0$. If
$\{f_k\}$ is a frame (Hilbert) for $\mathcal{H}_0$ space, then each
element of $\mathcal{H}_0$ can be recovered by an infinite
combinations of frame elements. That is, by the reconstruction
formula. If a signal $f$ is transmitted to a \break receiver, then
there are some kind of disturbances  in the received signal. To
\break overcome these disturbances from the receiver, frames plays
an important role. \break Actually, a signal in the space (after its
transmission) is in the form of the frame coefficients $\{<
f,S^{-1}f_k>\}$, $f \in \mathcal{H}_0$ . An error $\widetilde{e}$ is
always is expected with concern signal in the space. That is, actual
signal in the space is of the form $\{< f,S^{-1}f_k>+
\widetilde{e}\}$, where $\widetilde{e}$ is an error associated with
$f$. An interesting discussion in this direction is
given in  [13]. We extend the said problem to  Banach frames in general Banach spaces.\\

The following proposition provides sufficient condition for a
reconstruction system to satisfy property $\mathfrak{I}_\omega$ in
terms of non-existence of pre-frame operator associated with certain
error.

\begin{prop}
Suppose that $\{f^*_k\}$ has the reconstruction property for a
signal space$($Banach$)$ $\mathcal{X}$ with respect to $\{f_k\}$.
Let  $h^*$(error) be in $\mathcal{X}^*$
 for which there is no pre-frame operator
 $\mathfrak{S}_0$ such that  $(\cc{f^*_k+h^*}, \mathfrak{S}_0)$ is a Banach frame for
 $\mathcal{X}$, then $(\{f_k\}, \{f^*_k\})$ is a
 $\mathfrak{I}_\omega$-reconstruction system for $\mathcal{X}$.
\end{prop}

\proof Let $(\cc{f^*_k},\mathfrak{S})$ be an associated Banach frame
of  $(\{f_k\}, \{f^*_k\})$. If there exists no pre-frame operator
$\mathfrak{S}_0$ such that $(\cc{f^*_k+h^*}, \mathfrak{S}_0)$ is a
Banach frame for $\mathcal{X}$.
 Then,  there is a non-zero vector $f_0 \in \mathcal{X}$ such that $(f^*_k+ h^*)(f_0)=0$, for
all $k\in\NN$. By frame inequality of $(\cc{f^*_k},\mathfrak{S})$,
we conclude that  $h^*(f_0)\ne 0$. Put $\Psi=-\pi
(\dfrac{1}{h^*(f_0)}f_0)$.
 Then, $\Psi \in \mathcal{X}^{**}$ is such that $\Psi(f_k)=1$, for all $k\in\NN$.
Hence $(\{f_k\}, \{f^*_k\})$ is a
 $\mathfrak{I}_\omega$-reconstruction system for $\mathcal{X}$.
\endproof

\begin{rem}
The condition in Proposition $4.2$ is not necessary unless
 $\Psi$  correspond to a vector in $\mathcal{X}$. More precisely, we
 can find a certain error $h^* \in \mathcal{X}^*$ such that there exists
 no pre-frame operator $\mathfrak{S}_0$ associated with $\{f^*_k+
 h^*\}$ provided $\Psi \leftrightarrow f \in \mathcal{X}$.
\end{rem}

\begin{rem}
Let us continue with the outcomes in Proposition $4.2$, where
$(\{f_k\}, \{f^*_k\})$ is found to be a
$\mathfrak{I}_\omega$-reconstruction system for $\mathcal{X}$
provided there is no pre-frame operator $\mathfrak{S}_0$ such that
$(\cc{f^*_k+h^*}, \mathfrak{S}_0)$ is a Banach frame for
$\mathcal{X}$, where $h^*$ is certain choice of error(functional).
 A natural
problem arises, which is of determining a Banach space $\mathcal{B}$
for which the system $\{f^*_k+h^*\}$ admits a pre-frame operator.
Answer to this problem is positive provided $(\{f_k\}, \{f^*_k\})$
is pre-shrinking. The outline of construction of such a Banach space
can be understood as follows : Put $\zeta = -
\dfrac{1}{h^*(f_0)}f_0$ (where $f_0$ is same as in the proof of
Proposition $4.2$). Now, there is no pre-frame operator
$\mathfrak{S}_0$ associated with $\{f^*_k+h^*\}$, so there exists a
non-zero vector $g$ such that  $(f^*_k+ h^*)(g)=0$, for all
$k\in\NN$. By using frame inequality of the associated Banach frame
$(\cc{f^*_k},\mathfrak{S})$ we have $h^*(g) \ne 0$. Put
$\varphi=\dfrac{-1}{h^*(g)}g$. Then, $\varphi$ is a non-zero vector
in $\mathcal{X}$ such that $f^*_k(\varphi)=1$, for all $k\in\NN$.
Therefore, $f^*_k(\zeta - \varphi) = 0$ for all $k\in\NN$. Now
$(\{f_k\}, \{f^*_k\})$ is pre-shrinking, so we have $\zeta =
\varphi$. Hence $g = p \zeta$, where $p= -h^*(g)$. By using Lemma
$2.3$ there exists a pre-frame operator $\mathfrak{S}_1$ such that
$(\{f^*_k+ h^*\},\mathfrak{S}_1 )$ is a Banach frame(normalized
tight) for the Banach space $\mathcal{B}$ , where $\mathcal{B}^*=
[\zeta]^{\bot}$; $[\zeta]^{\bot} = \{f^* \in
\mathcal{X}^*:f^*(f)=0,f\in [\zeta]\}$.

\end{rem}
An application of Proposition $4.2$  is given below:

\begin{exa}
Let $(\{f_k\}, \{f^*_k\})$ be a  reconstruction system  given in
Example $3.4$ for $\mathcal{X}=c_0$. Then,
$\mathfrak{S}:\mathcal{X}_d=\{\{f^*_k(f)\}:f\in \mathcal{X}\}\to
\mathcal{X}$ is a bounded  linear operator such that
$(\cc{f^*_k},\mathfrak{S})$ is a  Banach frame $($associated$)$ for
$\mathcal{X}$ with respect to $\mathcal{X}_d$ and with bounds
$A=B=1$. Put $h^*=-f^*_4$ (this choice makes sense, because
disturbances  are not constant!). Then, $h^*$ is an error in
$\mathcal{X}^*$ for which there is no reconstruction operator
 $\mathfrak{S_0}$ such that
  $(\cc{f^*_k+ h^*, \mathfrak{S_0}})$ is a Banach frame for $\mathcal{X}$.
   Hence by Proposition $4.2$, $(\{f_k\}, \{f^*_k\})$ is a
   $\mathfrak{I}_\omega$-reconstruction system for $\mathcal{X}$. $\square$
\end{exa}

\begin{defn}
Fix $f \in \mathcal{X}$. A pair $(\{f_k\} ,\{f^*_k\})$, where
$(\{f_k\} \subset \mathcal{X},  \{f^*_k\} \subset \mathcal{X}^*)$ is
said to be \emph{localized at $f$}, if  $f =
\sum_{k=1}^{\infty}\epsilon_k f_k^*(f)f_k$, where $\{\epsilon_k\}$
is a sequence of scalars.
\end{defn}
 If $(\{f_k\} ,\{f^*_k\})$ is localized at
every $f \in \mathcal{X}$ with $\epsilon_k =1,$ for all $k$, then
$(\{f_k\} ,\{f^*_k\})$
 turns out to be a reconstruction system for $\mathcal{X}$.
 Consider a reconstruction system $(\{f_k\} ,\{f^*_k\})$ for $\mathcal{X}$ and
 $(\{f_k^*\}, \mathfrak{S}_0)$ be its associated Banach frame with respect to
 $\mathcal{X}_d$. Let $\mathrm{0} \ne \{\psi^*_k(f)\} \in \mathcal{X}_d, f \in \mathcal{X}$.  Then,
in  general, there is no pre-frame operator
$\widehat{\mathfrak{S}_0}$ associated with system
$\{\frac{1}{\psi^*_k(f)} f_k^* - \frac{1}{\psi^*_{k+1}(f)}
f_{k+1}^*\}$. This problem is also known as stability of
$(\{f_k^*\}, \mathfrak{S}_0)$ with respect to $\mathcal{X}_d$. If
$(\{f_k\} ,\{\psi^*_k\})$ is not localized at certain vectors in
$\mathcal{X}$, then we can find such pre-frame operator associated
with $\{\frac{1}{\psi^*_k(f)} f_k^* - \frac{1}{\psi^*_{k+1}(f)}
f_{k+1}^*\}$. This is what concluding proposition of this paper
says.
\begin{prop}
Let $(\{f_k\} ,\{f^*_k\})$ be a reconstruction system for
$\mathcal{X}$. Assume that $(\{f_k\} ,\{\psi^*_k\})$ is not
localized at $f \in \mathcal{X}$, where $\pi_f(\mathcal{E})=0;$
\break $\mathcal{E} = span(f^*_k)_k$. Then there exists  a pre-frame
operator $\widehat{\mathfrak{S}_0}$, such that \break
$(\{\frac{1}{\psi^*_k(f)} f_k^* - \frac{1}{\psi^*_{k+1}(f)}
f_{k+1}^*\}, \widehat{\mathfrak{S}_0})$ is a Banach frame for
$\mathcal{X}$.
\end{prop}

\proof Let $( \{f_k^*\}, \mathfrak{S})$ be associated Banach frame
of $(\{f_k\}, \{f^*_k\})$. Let, if possible, there is no
reconstruction operator  $ \widehat{\mathfrak{S}_0}$, such that
$(\{\frac{1}{\psi^*_k(f)} f_k^* - \frac{1}{\psi^*_{k+1}(f)}
f_{k+1}^*\},  \widehat{\mathfrak{S}_0})$ is a Banach frame for
$\mathcal{X}$. Then, there exists a non zero vector $f_0$ such that
\break  $(\frac{1}{\psi^*_k(f)} f_k^* - \frac{1}{\psi^*_{k+1}(f)}
f_{k+1}^*)(f_0) = 0,$ for all  $k  \in\mathbb{N}.$\\
 This gives
$$\frac{1}{\psi^*_k(f)}{f_k^*}(f_0)=
\frac{1}{\psi^*_{k+1}(f)}f_{k+1}^*(f_0),  \ \ k  \in\mathbb{N}.$$\\
By using frame inequality of $( \{f_k^*\},S)$,  we
obtain,$${f^*_k}(f_0)=\frac{\psi^*_k(f)}{\psi^*_1(f)}f^*_1(f_0)\neq0,
 \ \ k \in\mathbb{N}.$$\\
Since $(\{f_k\}, \{f_k^*\})$ is a  reconstruction system for
$\mathcal{X}$, we have
$$f_0=\sum_{k=1}^{\infty}f_k^*(f_0)f_k$$\\
$$=\sum_{k=1}^{\infty}\frac{\psi^*_k(f)}{\psi^*_1(f)}f^*_1(f_0)f_k.$$\\
Thus, $(\{f_k\} ,\{\psi^*_k\})$ is localized at $f_0 \in
\mathcal{X}$, where $\pi_{f_0}(\mathcal{E})=0$, a contradiction.
Hence there exists  a pre-frame operator $\widehat{\mathfrak{S}_0}$,
such that  $(\{\frac{1}{\psi^*_k(f)} f_k^* -
\frac{1}{\psi^*_{k+1}(f)} f_{k+1}^*\}, \widehat{\mathfrak{S}_0})$ is
a Banach frame for $\mathcal{X}$.
\endproof

\section{Conclusions}
 The notion of
 $\mathfrak{I}$-reconstruction property is proposed in section $3$ and
 its characterization in terms of frames in Banach spaces is given.
 More precisely, Proposition $3.5$ characterize $\mathfrak{I}$-reconstruction
property in terms  of existence of  pre-frame operator  but in a
contrapositive way. This situation is same as in electrodynamics,
where there is a game of movement of electron but charge given to
electron is negative! Moreover, the action of a functional from
$\mathcal{X}^{**}$ on a given system from $\mathcal{X}^*$  decide
the existence of pre-frame operator associated with certain system.
This looks like \emph{dynamics of reconstruction property}. By
motivation from the theory of frames for Hilbert spaces which
control the perturbed system associated   with a signal in
space(after its transmission), we extend the said situation to
Banach spaces. More precisely, Proposition $4.2$
 control the situation in abstract setting via non-existence of pre-frame operator.
  Finally, the notion of local  reconstruction system is proposed and its
  utility in complicated stability of associated Banach frames is reflected in Proposition $4.7$.

\mbox{}

\end{document}